\documentclass[11pt]{article}
\usepackage{amsmath, amssymb, amsthm, verbatim}

\begin{document}
\bibliographystyle{plain}

\newtheorem{theorem}{Theorem}
\newtheorem{lemma}[theorem]{Lemma}
\newtheorem{corollary}[theorem]{Corollary}

\newcommand\al{\alpha}
\newcommand\be{\beta}
\newcommand\ga{\gamma}
\newcommand\ff[1]{{\mathbb F_{#1}}} 

\title{Two characterizations of crooked functions}

\author{Chris Godsil\footnote{Department of Combinatorics and Optimization, University of Waterloo, Waterloo, ON N2L 3G1. email:~cgodsil@uwaterloo.ca. CG is supported by NSERC.} \and Aidan Roy\footnote{Institute for Quantum Information Science, University of Calgary, Calgary, AB, T2N 1N4. email:~aroy@qis.ucalgary.ca. AR is supported by NSERC and MITACS.}}
\maketitle\abstract{We give two characterizations of crooked functions: one based on the minimum distance of a Preparata-like code, and the other based on the distance-regularity of a crooked graph.}

\section{Introduction}

Highly nonlinear functions over finite vector spaces have attracted much interest in the last several years, for both their applications to cryptography (see~\cite{cha1} for example) and their connections to a variety of different combinatorial structures. The functions that are furthest from linear are called {\it perfect nonlinear}; unfortunately, none exist for binary vector spaces, which are the most cryptographically useful. However functions do exist in several lesser categories of nonlinearity, such as {\it almost perfect nonlinear}, {\it almost bent}, and {\it crooked}. We focus on the latter, which is the most specialized of the three. 

Crooked functions were introduced by Bending and Fon~Der~Flaass~\cite{ben1}, who, building on the graphs of de~Caen, Mathon and Moorhouse~\cite{dec2}, showed that every crooked function defines a distance-regular graph of diameter $3$ with a particular intersection array. Shortly thereafter, van~Dam and Fon~Der~Flaass~\cite{van3} observed that every crooked function defines a binary code of minimum distance $5$, similar to the classical Preparata code. In this paper, we show that the converse of each of these results is also true: crooked functions can be characterized using both Preparata-like codes (Theorem~\ref{thm:crookedcode}) and distance-regular graphs (Theorem~\ref{thm:crookedgraph}). Those codes and graphs offer a more combinatorial way of understanding the nature of nonlinear binary functions.


\section{Almost Perfect Nonlinear Functions}

Before considering crooked functions we need to characterize a more general class, namely {\it almost perfect nonlinear} functions. Throughout this article, let $V := V(m,2)$, a vector space of dimension $m$ over $\ff{2}$, with $m$ odd. Given a function $f:V \rightarrow V$, consider the following system of equations: 
\begin{equation}\label{eqn:apn}
\left\{ \begin{array}{rcl}
x + y & = & a \\
f(x) + f(y) & = & b
\end{array} \right\}.
\end{equation}
Note that solutions to \eqref{eqn:apn} come in pairs: if $(x,y)$ is a solution, then so is $(y,x)$. If $f$ is a linear function, then equation \eqref{eqn:apn} has $2^m$ solutions when $b = f(a)$. We say $f$ is {\it almost perfect nonlinear} if, for every $(a,b) \neq (0,0)$, the system has at most two solutions. Equivalently, $f$ is almost perfect nonlinear if and only if for all $a \neq 0$ in $V$, the set
\[
H_a(f) := \{f(x) + f(x+a) \mid x \in V\}
\]
has cardinality $2^{m-1}$. 

We may construct a binary code from a function on $V$ in the following manner. Identify $V$ with the finite field $\ff{2^m}$, and let $\al$ be a primitive element of $\ff{2^m}$. Also let $n = 2^m-1$, and assume $f:V \rightarrow V$ is a function such that $f(0) = 0$. We define a parity check matrix $H_f$ by 
\[
H_f := \left( \begin{matrix}
1 & \al & \al^2 & \ldots & \al^{n-1} \\
f(1) & f(\al) & f(\al^2) & \ldots & f(\al^{n-1}) 
\end{matrix} \right),
\]
and define the code $C_f$ to be the kernel of $H_f$ over $\ff{2}$. 

The code $C_f$ can be thought of as a generalization of the double error-correcting BCH code, which is the specific case of $f(\al) := \al^3$. It is clear from the parity check matrix that the minimum distance of $C_f$ is at least $3$, and it can be shown that the minimum distance is at most $5$. The following characterization is due to Carlet, Charpin, and Zinoviev~\cite[Theorem~5]{ccz}.

\begin{theorem}\label{thmparam}
The minimum distance of $C_f$ is $5$ if and only if $f$ is almost perfect nonlinear. In this case, the dimension of $C_f$ is
\[
k = 2^m - 2m - 1.
\] 
\end{theorem}

In the next section, we give a similar characterization of crooked functions, which are a special class of almost perfect nonlinear functions.

\section{Crooked Functions and Preparata-like Codes} 

A function $f: V \rightarrow V$ is {\it crooked} if the following three conditions hold: 
\begin{enumerate}
\item $f(0) = 0$; \label{itm:crooked1}
\item $f(x) + f(y) + f(z) \neq f(x+y+z) \quad$ for distinct $x$, $y$, and $z$; \label{itm:crooked2}
\item $f(x) + f(y) + f(z) \neq f(x+a) + f(y+a) + f(z+a) \quad$ for all $a \neq 0$. \label{itm:crooked3}
\end{enumerate}
Condition~\ref{itm:crooked2} is equivalent to almost perfect nonlinearity; thus every crooked function is almost perfect nonlinear. Condition~\ref{itm:crooked3} states that for every $a \neq 0$, no three points in $H_a(f)$ are collinear. It follows that $f$ is crooked if and only if $f(0) = 0$ and $H_a(f)$ is the complement of a hyperplane for all $a \neq 0$. Note that we are using the original definition of crooked functions given in~\cite{ben1}, rather than the generalization appearing in Byrne and McGuire~\cite{bm1} or Kyureghyan~\cite{kyu1}.

The canonical example of a crooked function is the {\it Gold} function. Identify $V$ with $\ff{2^m}$ for odd $m$; then $f(x) := x^{2^k+1}$ is called a Gold function if $\nobreak{\gcd(k,m) = 1}$. More generally, $f(x) := x^{2^k+2^j}$ is crooked provided that $\nobreak{\gcd(k-j,m) = 1}$, and Kyureghyan~\cite{kyu1} has shown that all crooked power functions have this form. For recent progress in constructing nonlinear functions which are not equivalent to the Gold functions, see~\cite{bcfl, bcp1, ekp1}. 


Just as almost perfect nonlinear functions give rise to BCH-like codes, crooked functions give to Preparata-like codes. Given $f: V \rightarrow V$ such that $f(0) = 0$, let $P_f$ be the code whose codewords are the characteristic vectors of $(S,T)$, for $S \subset V^*$ and $T \subset V$, such that the following three conditions hold:
\begin{itemize}
\item $|T|$ is even, 
\item ${\displaystyle \sum_{r \in S} r = \sum_{r \in T} r}$, and 
\item ${\displaystyle f\Big(\sum_{r \in S} r \Big) = \sum_{r \in S}f(r) + \sum_{r \in T}f(r)}$.
\end{itemize}
Identifying $V$ with $\ff{2^m}$, we get the actual Preparata code when $f(x)~:=~x^3$ and the generalized Preparata code when $f(x)~=~x^{2^k+1}$ (see~\cite{bak1}). In general $P_f$ is not linear, and it is easy to verify that $P_f$ always has minimum distance at least $3$. The following result is due to Van~Dam and Fon~Der~Flaass~\cite[Theorem~7]{van3}.

\begin{theorem}
If $f$ is crooked, then $P_f$ has minimum distance $5$ and size $2^{2^{m+1}-2m-2}$.
\label{thmcrkpf}
\end{theorem}

If $P_f$ has minimum distance $5$, then it is {\it nearly perfect}: it satisfies the Johnson bound~\cite[Theorem~17.13]{mac1} with equality. Hence $P_f$ has minimum distance at most $5$ for any $f$. We show the converse of Theorem~\ref{thmcrkpf}. 

\begin{theorem}
If $P_f$ has minimum distance $5$, then $f$ is crooked. 
\label{thm:crookedcode}
\end{theorem}

\begin{proof} 
We assumed in the definition of $P_f$ that $f(0) = 0$, so condition~\ref{itm:crooked1} of crookedness is satisfied. If $P_f$ has minimum distance $5$, then there is no pair $(\phi, T)$ in $P_f$ with $|T| = 4$. That is, for any distinct $w,x,y,z$ such that $w+x+y+z = 0$, 
\begin{equation}
f(w) + f(x) + f(y) + f(z) \neq 0.
\label{eqncrkapn}
\end{equation}
Thus condition~\ref{itm:crooked2} of crookedness is also satisfied, and it remains to show condition~\ref{itm:crooked3}. Since condition~\ref{itm:crooked2} is saitsfied, $f$ is almost perfect nonlinear and $C_f$ has dimension $2^m - 2m - 1$ by Theorem~\ref{thmparam}. But $C_f$ is the kernel of $H_f$, so it follows that $H_f$ has a column space of dimension $2m$, namely $V \times V$. This implies that for any $(a,b)$ in $V \times V$, there is a subset $S$ of $V^*$ such that 
\begin{equation}
\left(\sum_{r \in S} r, \sum_{r \in S} f(r)\right) = (a,b).
\label{eqncfspan}
\end{equation}

Given any $x \in V$, let $T = \{x,0\}$, so that $|T|$ is even and $\sum_{r \in T} r = x$. Then from equation \eqref{eqncfspan}, there exists some $S \subset V^*$ such that 
\[
\Big(\sum_{r \in S} r, \sum_{r \in S} f(r)\Big) = \left(x, 0 \right).
\] 
Choosing $S$ in this way, $(S,T)$ is in $P_f$. Now given any $y$, $z$ and $a \neq 0$, consider
\[
(S',T') := (S \oplus \{y\} \oplus \{y+a\}, T \oplus \{z\} \oplus \{z+a\}).
\]
This vector is at distance $4$ from $(S,T)$. Since $P_f$ has distance $5$, $(S',T')$ must not be in $P_f$. But $|T'|$ is even, and 
\[
\sum_{r \in S'} r = \sum_{r \in T'} r;
\]
hence for $(S',T') \notin P_f$ it must be the case that 
\[
f\Big(\sum_{r \in S'}r\Big) \neq \sum_{r \in S'}f(r) + \sum_{r \in T'} f(r).
\]
This implies
\[
f(x+a) \neq \sum_{r \in S} f(r) + f(y) + f(y+a) + \sum_{r \in T} f(r) + f(z) + f(z+a),
\]
or in other words
\[
f(x+a) \neq f(y) + f(y+a) + f(x) + f(z) + f(z+a).
\]
Thus condition~\ref{itm:crooked3} of crookedness is satisfied for $f$.
\end{proof}

\section{Crooked Graphs}

As usual, assume $f(0) = 0$. Define the {\it crooked graph} of $f$, denoted $G_f$, to have vertex set $V \times \ff{2} \times V$ with the following adjacency condition: distinct $(a,i,\al)$ and $(b,j,\be)$ are adjacent if and only if
\[
\al + \be = f(a + b) + (i+j+1)(f(a) + f(b)).
\]
It is not difficult to show that any two vertices in the subset 
\[
F_{ai} := \{ (a,i,\al) \mid \al \in V\}
\]
are at distance at least three, and that any two distinct subsets $F_{ai}$ and $F_{bj}$ are joined by a perfect matching. It follows that $G_f$ is a $2^m$-cover of the complete graph $K_{2^{m+1}}$, and each $F_{ai}$ is a fibre (for background on covers of complete graphs, see~\cite{god7}). The following theorem is given by Bending and Fon-Der-Flaass~\cite[Proposition 13]{ben1}.
\begin{theorem}
If $f$ is crooked, then $G_f$ is an antipodal distance-regular graph with intersection array
\[
\{ 2^{m+1} -1, 2^{m+1} -2, 1; 1, 2, 2^{m+1}-1 \}.
\]
\label{thmcrookedskew}
\end{theorem}

For background on distance-regular graphs, see~\cite{bcn1}. Again, we show the converse.

\begin{theorem}
If $G_f$ is distance-regular with intersection array
\[
\{ 2^{m+1} -1, 2^{m+1} -2, 1; 1, 2, 2^{m+1}-1 \},
\]
then $f$ is crooked.
\label{thm:crookedgraph}
\end{theorem}

\begin{proof}
For convenience, consider the graph $G'_f$ which consists of $G_f$ with a loop added to every vertex. This can be done by removing the restriction $(a,i,\al) \neq (b,j,\be)$ from the adjacency condition of $G_f$. If $G_f$ is distance-regular with $a_1 = 0$ and $c_2 = 2$, then $G'_f$ is a graph with the property  that any two vertices at distance $1$ or $2$ have exactly two common neighbours. 
That is, for any two vertices $(a,i,\al), (b,j,\be)$ such that $(a,i) \neq (b,j)$, there are exactly two vertices $(c,k,\ga)$ such that 
\begin{align}
\al + \ga & = f(a + c) + (i+k+1)(f(a) + f(c)), \label{eqn:cond1} \\
\be+\ga & = f(b + c) + (j+k+1)(f(b) + f(c)). \label{eqn:cond2}
\end{align}
We restrict our attention to the cases in which $i = j$, so that $a \neq b$. 
Adding \eqref{eqn:cond1} and \eqref{eqn:cond2} together, there are exactly two pairs $(c,k)$ such that
\[
\al + \be = f(a + c) + f(b + c) + (i+k+1)(f(a) + f(b)).
\]
Running over all values of $\al + \be$, we see that for fixed $(a,b,i)$, the multiset
\[
\{ f(a + c) + f(b + c) + (i+k+1)(f(a) + f(b)) \mid c \in V, k \in \ff{2} \}
\]
\begin{equation}
= \{ f(a+c) + f(b+c) \mid c \in V \} \cup \{ f(a+c) + f(b+c) + f(a) + f(b) \mid c \in V \}
\label{eqndd}
\end{equation}
contains each element of $V$ exactly twice. 

Now for some fixed $c$, consider $f(a+c) + f(b+c)$. Letting $c' := c + a + b$, we  have
\[
f(a+c) + f(b+c) = f(a+c') + f(b+c').
\]
However, the value $f(a+c) + f(b+c)$ only occurs twice in \eqref{eqndd}, so there is no third solution $c'' \neq c,c'$ such that
\[
f(a+c) + f(b+c) = f(a+c'') + f(b+c'').
\]
In other words, letting $x = a+c$, $y = b+c$, and $z = a+c''$, we have
\[
f(x) + f(y) \neq f(z) + f(x+y+z)
\]
for $z \neq x,y$. This is condition~\ref{itm:crooked2} of crookedness for $f$. Also because $\nobreak{f(a+c) + f(b+c)}$ has already occured twice in \eqref{eqndd}, there is no $c''$ such that 
\[
f(a+c) + f(b+c) = f(a+c'') + f(b+c'') + f(a) + f(b).
\]
Setting $x = a+c$, $y = a+c''$, $z = a$ and $w = a+b$, we have
\[
f(x) + f(x+w) \neq f(y) + f(y+w) + f(z) + f(z+w)
\]
for any $x,y,z$ and $w$, with $w \neq 0$. This is the condition~\ref{itm:crooked3} of crookedness, so $f$ is crooked.
\end{proof}

\bibliography{crooked}

\begin{thebibliography}{10}

\bibitem{bak1}
Ronald~D. Baker, Jacobus~H. van Lint, and Richard~M. Wilson.
\newblock On the {P}reparata and {G}oethals codes.
\newblock {\em IEEE Trans. Inform. Theory}, 29(3):342--345, 1983.

\bibitem{ben1}
T.~D. Bending and D.~Fon-Der-Flaass.
\newblock Crooked functions, bent functions, and distance regular graphs.
\newblock {\em Electron. J. Combin.}, 5(1):Research Paper 34, 14 pp.
  (electronic), 1998.

\bibitem{bcn1}
A.~E. Brouwer, A.~M. Cohen, and A.~Neumaier.
\newblock {\em Distance-Regular Graphs}.
\newblock Springer-Verlag, Berlin, 1989.

\bibitem{bcfl}
Lilya Budaghyan, Claude Carlet, Patrick Felke, and Gregor Leander.
\newblock An infinite class of quadratic apn functions which are not equivalent
  to power mappings.
\newblock In {\em Information Theory, 2006 IEEE International Symposium on},
  pages 2637--2641, 2006.

\bibitem{bcp1}
Lilya Budaghyan, Claude Carlet, and Alexander Pott.
\newblock New classes of almost bent and almost perfect nonlinear polynomials.
\newblock {\em IEEE Trans. Inform. Theory}, 52(3):1141--1152, 2006.

\bibitem{bm1}
Eimear Byrne and Gary McGuire.
\newblock On the non-existence of quadratic apn and crooked functions on finite
  fields.
\newblock 2005.
\newblock http://www.maths.may.ie/staff/gmg/APNniceWeilEBGMG.pdf.

\bibitem{ccz}
Claude Carlet, Pascale Charpin, and Victor Zinoviev.
\newblock Codes, bent functions and permutations suitable for {DES}-like
  cryptosystems.
\newblock {\em Des. Codes Cryptogr.}, 15(2):125--156, 1998.

\bibitem{cha1}
Florent Chabaud and Serge Vaudenay.
\newblock Links between differential and linear cryptanalysis.
\newblock In {\em Advances in cryptology---EUROCRYPT '94 (Perugia)}, volume 950
  of {\em Lecture Notes in Comput. Sci.}, pages 356--365. Springer, Berlin,
  1995.

\bibitem{dec2}
D.~de~Caen, R.~Mathon, and G.~E. Moorhouse.
\newblock A family of antipodal distance-regular graphs related to the
  classical {P}reparata codes.
\newblock {\em J. Algebraic Combin.}, 4(4):317--327, 1995.

\bibitem{ekp1}
Yves Edel, Gohar Kyureghyan, and Alexander Pott.
\newblock A new {APN} function which is not equivalent to a power mapping.
\newblock {\em IEEE Trans. Inform. Theory}, 52(2):744--747, 2006.

\bibitem{god7}
C.~D. Godsil and A.~D. Hensel.
\newblock Distance regular covers of the complete graph.
\newblock {\em J. Combin. Theory Ser. B}, 56(2):205--238, 1992.

\bibitem{kyu1}
Gohar Kyureghyan.
\newblock Crooked maps in finite fields.
\newblock In {\em 2005 {E}uropean {C}onference on {C}ombinatorics, {G}raph
  {T}heory and {A}pplications ({E}uro{C}omb '05)}, Discrete Mathematics \&
  Theoretical Computer Science Proceedings, AE, pages 167--170, 2005.

\bibitem{mac1}
F.~J. MacWilliams and N.~J.~A. Sloane.
\newblock {\em The Theory of Error-Correcting Codes}.
\newblock North-Holland Publishing Co., Amsterdam, 1977.

\bibitem{van3}
E.~R. van Dam and D.~Fon-Der-Flaass.
\newblock Uniformly packed codes and more distance regular graphs from crooked
  functions.
\newblock {\em J. Algebraic Combin.}, 12(2):115--121, 2000.

\end{thebibliography}

\end{document}